\documentclass[10pt,leqno,twoside]{article}

\usepackage[latin2]{inputenc}
\usepackage{amsmath}
\usepackage{amsfonts}
\usepackage{amssymb}
\usepackage{makeidx}

\newtheorem{statement}{Statement}
\newtheorem{theorem}{Theorem}

\newtheorem{proposition}{Proposition}

\newfont{\eurorm}{eurm10 scaled 1100}
\newfont{\eurorms}{eurm10 scaled 800}

\begin{document}

\bigskip
$\,$

\bigskip

\bigskip
\noindent
{\small \tt Preprint submitted to arXiv.org}

\vspace*{20.0 mm}

\begin{center}
{\large \bf
EXTENSIONS OF THE NATURAL APPROACH TO             \\[0.1 ex]
REFINEMENTS AND GENERALIZATIONS OF                \\[0.4 ex]
SOME TRIGONOMETRIC INEQUALITIES}
\end{center}

\smallskip

\begin{center}
{\em Branko Male\v sevi\' c${}^{\;\mbox{\scriptsize $\ast$}, \mbox{\tiny 1)}}\!$,
Tatjana Lutovac${}^{\,\mbox{\tiny 1)}}\!$,
Marija Ra\v sajski${}^{\,\mbox{\tiny 1)}}\!$,
Cristinel Mortici${}^{\,\mbox{\tiny 2)}}$}
\end{center}

\begin{center}
{\footnotesize \textit{${}^{1)}$School of Electrical Engineering,
University of Belgrade, \\[0.0ex]
Bulevar kralja Aleksandra 73, 11000 Belgrade, Serbia \\[2ex]
${}^{2)}$Valahia University of T\^{a}rgovi\c{s}te, Bd. Unirii 18, 130082 T\^{a}rgovi\c{s}te, Romania; \\[0.2ex]
Academy of Romanian Scientists, Splaiul Independen\c{t}ei 54, 050094 Bucharest, Romania;              \\[-0.2ex]
University Politehnica of Bucharest, Splaiul Independen\c{t}ei 313, 060042 Bucharest, Romania}}
\end{center}

\medskip
\noindent {\small \textbf{Abstract.}
In this paper we propose a new method for sharpening and refinements of some trigonometric inequalities.
We apply these ideas to some inequalities of {\sc Wilker}-{\sc Cusa}-{\sc Huygens}'s type.



\footnote{$\!\!\!\!\!\!\!\!\!\!\!\!\!\!$
{\scriptsize
${}^{\mbox{\scriptsize $\ast$}}$Corresponding author.             \\[0.0 ex]
Emails:                                                           \\[0.0 ex]
{\em Branko Male\v sevi\' c} {\tt $<$branko.malesevic@etf.rs$>$},
{\em Tatjana Lutovac} {\tt $<$tatjana.lutovac@etf.rs$>$},         \\[0.0 ex]
{\em Marija Ra\v sajski} {\tt $<$marija.rasajski@etf.rs$>$},
{\em Cristinel Mortici} {\tt $<$cristinel.mortici@valahia.ro$>$}}}

\vspace*{-2.5 mm}

{\small \tt MSC: Primary 33B10; Secondary 26D05}

\section{Introduction}

Inequalities involving trigonometric functions are used in many applications in various fields of mathematics such as
difference equations and inequalities \cite{R_P_Agarwal_1992}, theory of stability, theory of approximations, etc.
A method called {\it the natural approach}, proposed by {\sc Mortici} in \cite{Mortici_2011}, uses the idea of comparing functions to their corresponding {\sc Taylor} polynomials. This method has been successfully applied to prove and approximate
a wide category of trigonometric inequalities \cite{G_V_Milovanovic_2014}.

In this paper we extend the ideas of the natural approach
by comparing and replacing functions with  their corresponding power series.
In particular, we focus on the results of {\sc Mortici} in
\cite{Mortici_2011} related to {\sc Wilker}-{\sc Cusa}-{\sc Huygens}'s inequalities and
give generalizations and refinements of the inequalities stated in
Theorems 1, 2, 3, 4, 5 and 6 in that paper. They are cited below for readers  convenience.

\begin{statement}
\label{statement_1}
$(${\rm \cite{Mortici_2011}}, Theorem $1$$)$
For every $0<x<\pi/2$, we have
$$
\displaystyle - \frac{x^{4}}{15}\, < \,
\cos{x} - \left( \frac{\sin{x}}{x} \right)^{\!3}\!
\, < \,
- \frac{x^{4}}{15} + \frac{23x^{6}}{1890}.
$$
\end{statement}

\begin{statement}
\label{statement_2}
$(${\rm \cite{Mortici_2011}},  Theorem $2$$)$
For every $0<x<\pi/2$, we have
$$
\displaystyle - \frac{1}{180}x^{4}\, < \,
\frac{\sin{x}}{x} -\frac{ \cos{x} + 2}{3}\, < \,
- \frac{1}{180}x^{4} + \frac{1}{3780}x^{6}.
$$
\end{statement}

\begin{statement}
\label{statement_3}
$(${\rm \cite{Mortici_2011}}, Theorem $3$$)$
For every $0<x<\pi/2$, we have
$$
\displaystyle
3 + \left(\frac{3}{20}x^{4} - \frac{3}{140}x^{6} \right)\frac{1}{\cos{x}}\, < \,
2 \frac{\sin{x}}{x} \, + \, \frac{\mbox{\rm tan}\,x}{x}\, < \,
3 \, + \, \frac{3}{20}\, \frac{x^{4}}{\cos{x}}.
$$
\end{statement}

\begin{statement}
\label{statement_4}
$(${\rm \cite{Mortici_2011}}, Theorem $4$$)$ For every $0<x<\pi/2$, we have
$$
\displaystyle 2 + \left(\frac{8 \,x^{4}}{45}- \frac{8 \,x^{6}}{105}\right)\frac{1}{\cos{x}}\, < \,
\left(\frac{\sin{\,x}}{x}\right)^{\!2}\, + \, \frac{\tan{x}}{x}\, < \,
2 + \frac{8 \,x^{4}}{45}\frac{1}{\cos{x}}
.
$$
\end{statement}

\begin{statement}
\label{statement_5} $(${\rm \cite{Mortici_2011}}, Theorem $5$$)$ For every $0<x<\pi/2$, we have
$$
\displaystyle
\left(\frac{x}{\sin{x}}\right)^{\!2} + \frac{x}{\tan{x}} > 2 + \frac{2}{45}\,x^4.
$$
\end{statement}

\begin{statement}
\label{statement_6} $(${\rm \cite{Mortici_2011}}, Theorem $6$$)$ For every $0<x<\pi/2$, we have
$$\displaystyle
3 \frac{x}{\sin{x}} + \cos{x} > 4 + \frac{1}{10}\,x^4 + \frac{1}{210}x^6.
$$
\end{statement}

\medskip
\section{Preliminaries}


First, let us recall some of the well-known power series expansions that will be used in our proofs.

\smallskip
For $x \in R$ the following power series expansions hold:
\begin{equation}
\label{sin}
\mbox{\rm sin}\,x \,= \sum\limits_{k=0}^{\infty}{(-1)^{k}\displaystyle\frac{1}{(2k+1)!}\, x^{2k+1}},
\quad \quad
\mbox{\rm cos}\, x \,= \sum_{k=0}^{\infty}(-1)^{k}{\displaystyle\frac{1}{ \, (2\,k)!}\, x^{2k}}.
\end{equation}

Also, according to \cite{Gradshteyn-Ryzhik}, for $x \in R$ we have the following power series expansions:
\begin{equation}
\label{cos3} \mbox{\rm cos}^{3}\, x \,= \, \frac{1}{4}\, \sum_{k=1}^{\infty}{(-1)^{k}
 \displaystyle\frac{3^{2k} + 3 }{ \, (2\,k)!}\, x^{2k}},
\end{equation}
and
\begin{equation}
\label{sin3}
\sin^3 x
=
\frac{1}{4}
\sum_{k=1}^{\infty}{(-1)^{k+1}\displaystyle\frac{3^{2k+1}-3}{(2k+1)!}x^{2k+1}}.
\end{equation}

For $\displaystyle x \!\in\! \left( 0, \displaystyle\frac{\pi}{2} \right)$, according to \cite{Gradshteyn-Ryzhik},
the following series expansions hold:
\begin{equation}
\label{cosec}
{\mathop{\rm cosec}\nolimits} (x) = \frac{1}{x} + \mathop \sum \limits_{k = 1}^{  \infty } \frac{{\left|
{{\mbox{\boldmath $B$}_{2k}}} \right|({2^{2k}} - 2)}}{{(2k)!}}{x^{2k - 1}},
\end{equation}
\begin{equation}
\label{cosec2}
{{\mathop{\rm cosec}\nolimits} ^2}(x) = \frac{1}{{{x^2}}} + \mathop \sum \limits_{k = 1}^{  \infty }
\frac{{\left| {{\mbox{\boldmath $B$}_{2k}}} \right|(2k - 1){4^k}}}{{(2k)!}}{x^{2k - 2}},
\end{equation}
and
\begin{equation}
\label{cotan} {\mathop{\rm cotan}\nolimits} (x) = \frac{1}{x} -
\mathop \sum \limits_{k = 1}^{ \infty } \frac{{\left| {{\mbox{\boldmath $B$}_{2k}}}
\right|{4^k}}}{{(2k)!}}{x^{2k - 1}}
\end{equation}
where $\mbox{\boldmath $B$}_{i}$ are {\sc Bernoulli}'s numbers.

\bigskip
\noindent
{\bf Theorem WD.}  $(\,${\rm \cite{Wu_Debnath_2009}}, Theorem $2$$)$

\smallskip
\noindent
{\em Suppose that $f(x)$ is a real function on $(a,b)$, and that $n$ is a positive integer such that $f ^{(k)}(a+),
f^{(k)}(b-)$,
$\left(k \!\in\! \{0,1,2, \ldots ,n\}\right)$ exist.

\medskip
\noindent
{\boldmath $(i)$} Supposing that  $(-1)^{(n)} f^{(n)}(x)$~is~in\-cre\-asing on $(a,b)$, then
for all $x \in (a,b)$ the following inequality holds$:$
\begin{equation}
\label{Debnath_Wu_first}
\begin{array}{c}
\displaystyle\sum_{k=0}^{n-1}{\mbox{\small $\displaystyle\frac{f^{(k)}(b\mbox{\footnotesize
$-$})}{k!}$}(x-b)^k}
+
\frac{1}{(a-b)^n}
{\bigg (}\!
f(a\mbox{\footnotesize $+$})
-
\displaystyle\sum_{k=0}^{n-1}{\mbox{\small $\displaystyle\frac{(a-b)^{k}f^{(k)}(b\mbox{\footnotesize
$-$})}{k!}$}
\!{\bigg )} (x-b)^{n}}
\\[2.0 ex]
<
f(x)
<
\displaystyle\sum_{k=0}^{n}{\frac{f^{(k)}(b\mbox{\footnotesize $-$})}{k!}(x-b)^{k}}.
\end{array}
\end{equation}
Furthermore, if $(-1)^{n} f^{(n)}(x)$ is decreasing on $(a,b)$, then the reversed inequality of {\rm
(\ref{Debnath_Wu_first})} holds.

\medskip
\noindent
{\boldmath $(ii)$} Supposing that  $f^{(n)}(x)$ is increasing on $(a,b)$, then for all $x \!\in\! (a,b)$ the following
inequality  holds:
\begin{equation}
\label{Debnath_Wu_second}
\begin{array}{c}
\displaystyle\sum_{k=0}^{n-1}{\mbox{\small $\displaystyle\frac{f^{(k)}(a\mbox{\footnotesize
$+$})}{k!}$}(x-a)^k}
+
\frac{1}{(b-a)^n}
{\bigg (}\!
f(b-)
-
\displaystyle\sum_{k=0}^{n-1}{\mbox{\small $\displaystyle\frac{(b-a)^{k}f^{(k)}(a\mbox{\footnotesize
$+$})}{k!}$}
\!{\bigg )} (x-a)^{n}}
\\[2.0 ex]
>
f(x)
>
\displaystyle\sum_{k=0}^{n}{\frac{f^{(k)}(a\mbox{\footnotesize $+$})}{k!}(x-a)^{k}}.
\end{array}
\end{equation}
Furthermore, if $f^{(n)}(x)$ is decreasing on $(a,b)$, then the reversed inequality~of~\mbox{\rm
(\ref{Debnath_Wu_second})} holds.}

\smallskip
\noindent
Let us mention that an interesting application of Theorem WD is~given~in~\cite{MRL_2017c},
see also~\cite{Milica_Makragic_2017}.

\break

\section{Main results}

\medskip
 We need the following Theorem  for the proofs of Theorems 1, 2, 3 and 4.

\smallskip
\begin{proposition}\label{Leibniz} ~\\
Let the series $\displaystyle f(x)= \sum_{k=1}^{\infty}(-1)^{k} A(k) x^{2k}$ converges for
$x \!\in\! (0,c), \,\, c \!\in\! R^{+}\!.$ Suppose the following statements are true$\,:$
\begin{itemize}
\item[$\,(i)\,$] If $c\!<\!1$ then the sequence $\left\{ A(k) \right\}_{k \in N}$ is a positive decreasing sequence
that converges to $0$.
\item[$(ii)$] If $c \!\geq\! 1$ then the sequence $\displaystyle \left\{ A(k) \right\}_{k \in N}$ is a positive sequence, \mbox{$\!\!\lim\limits_{k\rightarrow +\infty}\!\!c^kA(k)\!=\!0$} and $\displaystyle A(k) \!>\! c^2 A(k+1)$ for $k \!\geq\! 1$.
\end{itemize}

\noindent
Then for all $x \!\in\! (0, c)$ and for all $n \!\in\! N$ and $m \!\in\! N$ we have:
\begin{equation}
\label{Leibniz_03}
\displaystyle\sum_{k=1}^{2n}{(-1)^{k+1} A(k) \, x^{2k}}
<
f(x)
<
\displaystyle\sum_{k=1}^{2n+1}{(-1)^{k+1} A(k) \, x^{2k}}
\end{equation}
and
\begin{equation}
\label{Leibniz_04}
{\Big |}\,f(x) - \displaystyle\sum_{k=1}^{m}{(-1)^{k} A(k) \, x^{2k}}\,{\Big |}
<
A(m+1) \, x^{2m+2}
<
c^{2m+2}\, A(m+1).
\end{equation}
\end{proposition}
{\bf Proof.}
Suppose that $c \!<\! 1$. Then for every $x \!\in\! (0,c)$ the positive
sequence $\left\{A(k) \, x^{2k} \right\}_{k \in N}$ decreases monotonically and
\mbox{$\displaystyle \lim _{k\rightarrow \infty } A(k)\, x^{2k} =0$.} Thus,  assertions
(\ref{Leibniz_03}) and (\ref{Leibniz_04}) immediately follow from {\sc Leibniz}'s theorem for
the alternating series.

\smallskip
Suppose now that $c \!\geq\! 1$. We have
$$
\displaystyle
\begin{array}{lcl}\displaystyle
f(x) \!&\!\!=\!\!&\! \displaystyle \sum_{k=1}^{\infty}(-1)^{k+1} A(k) x^{2k}     \\[1.5 ex]
 \!&\!\!=\!\!&\! \displaystyle \sum_{k=1}^{\infty}(-1)^{k+1} A(k) \, c^{2k}\, \left(\frac{x}{c}\right)^{\!2k}.
\end{array}
$$
Let us introduce the substitution $t\!=\!\displaystyle\frac{x}{c}$ in the previous power series and consider the series
\begin{equation}
 \label{smena-t}
 \displaystyle
 \sum_{k=1}^{\infty}(-1)^{k+1} A(k) c^{2k} t^{2k} , \quad \mbox{for}\,\,\, t \!\in\! (0,1).
 \end{equation}
For the assumption $\displaystyle A(k) > c^2A(k+1)$ we have
$$
\displaystyle
A(k) > c^2 A(k+1)
\,\,\, \Longleftrightarrow \,\,\,
A(k)c^{2k} > c^{2k+2} A(k+1),
~~ \mbox{for} ~ k \!\geq\! 1.
$$
Hence, we conclude that for every $t \!\in\! (0,1)$  power series (\ref{smena-t}) satisfies {\sc Leibniz}'s
theorem for the alternating series, and for all $n, m \in N$ we have:
\begin{equation}
\label{1}
{\small
\displaystyle\sum_{k=1}^{2n}{(-1)^{k+1}\, A(k) \, c^{2k} \, t^{2k}}
<
\sum_{k=1}^{\infty}(-1)^{k+1}\, A(k) \,c^{2k}\, t^{2k}
<
\displaystyle\sum_{k=1}^{2n+1}{(-1)^{k+1} \, A(k) \, c^{2k} \, t^{2k}}
}
\end{equation}
and
\begin{equation}
\label{2}
\begin{array}{lcl}
{\Big |}\,
\displaystyle \sum_{k=1}^{\infty}(-1)^{k+1} A(k) c^{2k} t^{2k}
-
\displaystyle\sum_{k=1}^{m}{(-1)^{k+1}\, A(k) \, c^{2k}\, t^{2k}}\,{\Big |}
\!&\!\!<\!\!&\! A(m+1) \, c^{2m+2}\, t^{2m+2}   \\ [1.25 ex]
\!&\!\!<\!\!&\! A(m+1) \, c^{2m+2}.
\end{array}
\end{equation}
Returning the variable $x= t\, c$ to (\ref{1}) and (\ref{2}) gives the assertions of proposition.
\hfill $\Box$

\bigskip


\medskip
\subsection{ Refinements of the inequalities in Statement \mbox{1}}

\smallskip
We propose the following improvement and generalization of Statement 1:

\smallskip
\begin{theorem} ~\\
\label{improvement-s-1}
$(i)$ For every $\displaystyle \, x \!\in\! \left(0, \mbox{\small $\displaystyle\frac{\pi}{2}$}\right)$ and every $n \in N$, we have
\begin{equation}
\label{t-1}
\begin{array}{c}
\displaystyle
\sum_{k=2}^{2n} (-1)^{k} A(k) \, x^{2k}
<
\cos x - \left(\displaystyle\frac{\sin x}{x}\right)^{\!3}
<
\displaystyle
\sum_{k=2}^{2n+1}(-1)^{k} A(k) \, x^{2k},
\end{array}
\end{equation}
where
$\displaystyle
A(k) = \mbox{\small $\displaystyle \frac{3^{2k+3}\!-\!32k^3\!-\!96k^2\!-\!88k\!-\!27}{4\,(2k+3)!}$}.
$

\smallskip
\noindent $(ii)$ For every $x \!\in\! \left(0, \mbox{\small $\displaystyle\frac{\pi}{2}$}\right)$ and every $m \in N$,
we have the following error estimation:
\begin{equation}
\begin{array}{l}
{\Big |}\,
\cos x
-
\left(\displaystyle\frac{\sin x}{x}\right)^{\!3}
\!-
\displaystyle\sum_{k=1}^{m}{(-1)^{k}
A(k) \, x^{2k}}\,{\Big |}\,
<\, A(m+1) \, x^{2m+2}.
\end{array}
\end{equation}
\end{theorem}

\noindent
{\bf Examples. }   \\
Let $x \!\in\! \left(0, \mbox{\small $\displaystyle\frac{\pi}{2}$}\right)$. \\
For $n=1$ we get Statement~\ref{statement_1}. \\[0.5 ex]
For $n>1$ we have the following new results:

\begin{itemize}
\item Taking $n=2$ in (\ref{t-1}) gives:

$
- \mbox{\small $\displaystyle\frac{1}{15}$}x^{4}
+ \mbox{\small $\displaystyle\frac{23}{1890}$}x^{6}
- \mbox{\small $\displaystyle\frac{41}{37800}$}x^{8}\, < \,
\cos x \, - \left(\displaystyle\frac{\sin x}{x}\right)^{\!3} <
$                                                                                                             \\

\hfill
$
\, < \,
- \mbox{\small $\displaystyle\frac{1}{15}$}x^{4}
+ \mbox{\small $\displaystyle\frac{23}{1890}$}x^{6}
- \mbox{\small $\displaystyle\frac{41}{37800}$}x^{8}
+ \mbox{\small $\displaystyle\frac{53}{831600}$}x^{10}
$

\item Taking $n=3$ in (\ref{t-1}) gives:

$
\begin{array}{l}
\displaystyle
- \mbox{\small $\displaystyle\frac{1}{15}$}x^{4}
+ \mbox{\small $\displaystyle\frac{23}{1890}$}x^{6}
- \mbox{\small $\displaystyle\frac{41}{37800}$}x^{8}
+ \mbox{\small $\displaystyle\frac{53}{831600}$}x^{10}
- \mbox{\small $\displaystyle\frac{74677}{27243216000}$}x^{12} <                                              \\[0.5 em]
< \, \cos x - \left(\mbox{\small $\displaystyle\frac{\sin x}{x}$}\right)^{\!3} <                              \\[0.5 em]
<
-
\mbox{\small $\displaystyle\frac{1}{15}$}x^{4}
+
\mbox{\small $\displaystyle\frac{23}{1890}$}x^{6}
-
\mbox{\small $\displaystyle\frac{41}{37800}$}x^{8}
+
\mbox{\small $\displaystyle\frac{53}{831600}$}x^{10}
-
\mbox{\small $\displaystyle\frac{74677}{27243216000}$}x^{12}
+
\mbox{\small $\displaystyle\frac{989}{10897286400}$}x^{14}
\end{array}$

\item Taking $n=4$ in (\ref{t-1}) gives:

$
-
\mbox{\small $\displaystyle\frac{x^{4}}{15}$}
+
\mbox{\small $\displaystyle\frac{23 \,x^{6}}{1890}$}
-
\mbox{\small $\displaystyle\frac{41\,x^{8}}{37800}$}
+
\mbox{\small $\displaystyle\frac{53\,x^{10}}{831600}$}
-
\mbox{\small $\displaystyle\frac{74677\,x^{12}}{27243216000}$}
+
\mbox{\small $\displaystyle\frac{989\,x^{14}}{10897286400}$}
-
\mbox{\small $\displaystyle\frac{79649\,x^{16}}{33345696384000}$}
\, < \,
$ \\[0.5 em]
$
< \,
\cos x - \left(\mbox{\small $\displaystyle\frac{\sin x}{x}$}\right)^{\!3}
\, < \,
-
\mbox{\small $\displaystyle\frac{x^{4}}{15}$}
+
\mbox{\small $\displaystyle\frac{23\,x^{6}}{1890}$}
-
\mbox{\small $\displaystyle\frac{41\,x^{8}}{37800}$}
+
\mbox{\small $\displaystyle\frac{53\,x^{10}}{831600}$}
-
\mbox{\small $\displaystyle\frac{74677\,x^{12}}{27243216000}$}
\, +
$ \\ [0.5 em]
$
+ \,
\mbox{\small $\displaystyle\frac{989\,x^{14}}{10897286400}$}
-
\mbox{\small $\displaystyle\frac{79649\,x^{16}}{33345696384000}$}
+
\mbox{\small $\displaystyle\frac{454007\,x^{18}}{8869955238144000}$}
$

\item[\mbox{etc.}]
\end{itemize}

\bigskip
\noindent
{\bf Proof of Theorem~\ref{improvement-s-1}.}

\smallskip
Consider the function
$$
f(x) = \cos x - \left(\displaystyle\frac{\sin x}{x}\right)^{\!3},
\quad \mbox{for}\,\, x \!\in\! \left(0, \mbox{\small $\displaystyle\frac{\pi}{2}$}\right).
$$
Based on  power series expansion (\ref{sin}) and (\ref{sin3}) we have
$$
f(x)
=
\displaystyle\sum_{k=1}^{\infty}{(-1)^{k} A(k) \, x^{2k}}
=
\displaystyle\sum_{k=2}^{\infty}{(-1)^{k} A(k) \, x^{2k}},
$$
for all $x \!\in\! R$, where$:$
\begin{equation}
\label{Def_a_k}
\begin{array}{rcl}
A(k)
\!&\!\!=\!\!&\!
\displaystyle\frac{1}{4}\frac{3^{2k+3}-3}{(2k+3)!} - \displaystyle\frac{1}{(2k)!}  \\[2.1 ex]
\!&\!\!=\!\!&\!
\displaystyle\frac{3^{2k+3}-32k^3-96k^2-88k-27}{4(2k+3)!}.
\end{array}
\end{equation}
For $ c=\pi/2$ we have:
 $$ A(k)>0 \;\, \mbox{for} \; k\geq 2 \quad
 \mbox{and} \quad \displaystyle \lim _{k\rightarrow \infty }c^{2k}A(k) =0.$$
Also,
$$
\begin{array}{l}\displaystyle
c^{2}\,\mathrm{A}(k + 1) < \mathrm{A}(k) ~ \Longleftrightarrow ~  \\[0.25 em]
32\,k^{5} \!+\! 240\,k^{4} \!+\! (680 \!-\! 8\,c^{2})\,k^{3} \!+\!
(900 \!-\! 45\,c^{2} )\,k^{2} \!+\! (548 \!-\! {\displaystyle \frac{161\,c^{2}}{2}} )
\,k \!-\! 45\,c^{2} \!+\! {\displaystyle \frac{477}{4}}   \\[1.0 em]
\hspace{35mm} < \, {\displaystyle \frac{3\, (4\,c^{2}\,k^{2} + 18\,c^{2}\,k + 20\,c
^{2} - 9)\,}{4}\,9^{k+1}}.
\end{array}
$$
As the last inequality holds for $k \geq 1$, the assertions of Theorem~\ref{improvement-s-1} immediately follows
from Proposition~\ref{Leibniz}. \hfill $\Box$


\medskip
\subsection{ Refinements of the inequalities in Statement \mbox{2}}

\smallskip
We propose the following improvement and generalization of Statement 2:

\smallskip
\begin{theorem} ~\\
\label{improvement-s-2}
 $(i)$ For every $\displaystyle \, x \!\in\! \left(0, \mbox{\small $\displaystyle\frac{\pi}{2}$}\right)$ and every $n \in N$, we have
\begin{equation}
\label{t-2}
\begin{array}{c}
\displaystyle
\sum_{k=2}^{2n} (-1)^{k+1} B(k) \, x^{2k}
<
 \frac{\sin{x}}{x} -\frac{ \cos{x} + 2}{3}
<
\displaystyle
\sum_{k=2}^{2n+1}(-1)^{k+1} B(k) \, x^{2k},
\end{array}
\end{equation}
where
$\displaystyle
B(k) = {\displaystyle \frac{2}{3}}\,{\displaystyle \frac{k - 1}{(2
\,k + 1)\mathrm{!}}}.
$

\smallskip
\noindent
$\,(ii)$ For every $x \!\in\! \left(0, \mbox{\small $\displaystyle\frac{\pi}{2}$}\right)$ and every $m \in N$,
we have the following error estimation:
\begin{equation}
\begin{array}{l}
{\Big |}\,\displaystyle
\frac{\sin{x}}{x} -\frac{ \cos{x} + 2}{3}
\!-\!
\displaystyle\sum_{k=0}^{m}{(-1)^{k+1}
B(k) \, x^{2k}}\,{\Big |} \, < \,B(m+1) \, x^{2m+2}.
\end{array}
\end{equation}
\end{theorem}
{\bf Examples.} \\
Let $ x \!\in\! \left(0, \mbox{\small $\displaystyle\frac{\pi}{2}$}\right)$. \\
For $n = 1$ we get Statement~\ref{statement_2}. \\[1.0 ex]
For $n >1 $ we have the following new results:
\begin{itemize}
\item Taking $n=2$ in (\ref{t-2}) gives:

$
-
\mbox{\small $\displaystyle\frac{1}{180}$}\,x^{4}
+
\mbox{\small $\displaystyle\frac{1}{3780}$}\,x^{6}
-
\mbox{\small $\displaystyle\frac{1}{181440}$}\,x^{8}
<
\mbox{\small $\displaystyle\frac{\sin x}{x}$}
-
\mbox{\small $\displaystyle\frac{1}{3}$}\cos x
-
\mbox{\small $\displaystyle\frac{2}{3}$} <
$
$
-
\mbox{\small $\displaystyle\frac{1}{180}$}\,x^{4}
+
\mbox{\small $\displaystyle\frac{1}{3780}$}\,x^{6}
-
-\mbox{\small $\displaystyle\frac{1}{181440}$}\,x^{8}
+
\mbox{\small $\displaystyle\frac{1}{14968800}$}\,x^{10}
$

\item Taking $n=3$ in (\ref{t-2}) gives:

$
-
\mbox{\small $\displaystyle\frac{\,x^{4}}{180}$}
+
\mbox{\small $\displaystyle\frac{\,x^{6}}{3780}$}
-
\mbox{\small $\displaystyle\frac{\,x^{6}}{181440}$}
+
\mbox{\small $\displaystyle\frac{\,x^{10}}{14968800}$}
-
\mbox{\small $\displaystyle\frac{\,x^{12}}{1868106240}$}
< \!
\mbox{\small $\displaystyle\frac{\sin \,x}{x}$}\!
-\!
\mbox{\small $\displaystyle\frac{1}{3}$}\cos x
-\!
\mbox{\small $\displaystyle\frac{2}{3}$}\! <
$
$
-
\mbox{\small $\displaystyle\frac{\,x^{4}}{180}$}
+
\mbox{\small $\displaystyle\frac{\,x^{6}}{3780}$}
-
\mbox{\small $\displaystyle\frac{\,x^{8}}{181440}$}
+
\mbox{\small $\displaystyle\frac{\,x^{10}}{14968800}$}
-
\mbox{\small $\displaystyle\frac{\,x^{12}}{1868106240}$}
+
\mbox{\small $\displaystyle\frac{\,x^{14}}{326918592000}$}
$

\item Taking $n=4$ in (\ref{t-2}) gives:

$
-
\mbox{\small $\displaystyle\frac{1}{180}$}\,x^{4}
+
\mbox{\small $\displaystyle\frac{1}{3780}$}\,x^{6}
-
\mbox{\small $\displaystyle\frac{1}{181440}$}\,x^{8}
+
\mbox{\small $\displaystyle\frac{1}{14968800}$}\,x^{10}
-
\mbox{\small $\displaystyle\frac{1}{1868106240}$}\,x^{12}
+
\mbox{\small $\displaystyle\frac{1}{326918592000}$}\,x^{14} -
-\mbox{\small $\displaystyle\frac{1}{76218734592000}$}\,x^{16} <
\mbox{\small $\displaystyle\frac{\sin x}{x}$}
-
\mbox{\small $\displaystyle\frac{1}{3}$}\cos x
-
\mbox{\small $\displaystyle \frac{2}{3}$}
<
-
\mbox{\small $\displaystyle\frac{1}{180}$}\,x^{4}
+
\mbox{\small $\displaystyle\frac{1}{3780}$}\,x^{6}
-
\mbox{\small $\displaystyle\frac{1}{181440}$}\,x^{8}
+
+\mbox{\small $\displaystyle\frac{1}{14968800}$}\,x^{10}
-
\mbox{\small $\displaystyle\frac{1}{1868106240}$}\,x^{12}
+
\mbox{\small $\displaystyle\frac{1}{326918592000}$}\,x^{14}
-
\mbox{\small $\displaystyle\frac{1}{76218734592000}$}\,x^{16}
+
+\mbox{\small $\displaystyle\frac{1}{22808456326656000}$}\,x^{18}
$
\item[etc.]
\end{itemize}

\noindent
{\bf Proof of Theorem~\ref{improvement-s-2}.}

\smallskip
Consider the function $$f(x)= {\displaystyle \frac{\mathrm{sin}\,x
}{x}} - {\displaystyle \frac{1}{3}}\,\mathrm{cos}\,x -
{\displaystyle \frac{2}{3}}, \quad \mbox{for} \; x \in \left(0, \mbox{\small $\displaystyle\frac{\pi}{2}$}\right).$$

\noindent
Based on power series expansion (\ref{sin}) we have
$$
f(x)
= \displaystyle
- {\displaystyle \frac{2}{3}
} + {\displaystyle \sum _{k=0}^{\infty}}\,(-1)^{k + 1}\,  {\displaystyle \frac{2}{3}}\,{\displaystyle \frac{k - 1}{(2 \,k + 1)\mathrm{!}}} \,x^{2\,k}\, = \, {\displaystyle \sum _{k=2}^{\infty}}\,(-1)^{k + 1}\,  {\displaystyle \frac{2}{3}}\,{\displaystyle \frac{k - 1}{(2 \,k + 1)\mathrm{!}}} \,x^{2\,k}.
$$

The sequence $\left\{ B(k) \right\}_{k\in N, \,k \geq 2}$
satisfies the recurrence relation:
$$
B(k + 1)= {\displaystyle
\frac{k}{2\, (k^{2} - 1)\,(2\,k + 3)}}\, B(k).
$$
For $c=\pi/2$ we have:
$$
B(k)>0 \;\; \mbox{for} \;\; k\geq 2, \; \;
\quad\mbox{and} \;\; \quad \displaystyle \lim _{k\rightarrow \infty }c^{2k}B(k) =0.
$$
Also,
$$
\begin{array}{lcl}\displaystyle
c^{2}\,B(k + 1) < B(k)
\!&\!\!\Longleftrightarrow\!\!&\!
\left( {\displaystyle
\frac{c^{2}\, k}{2\, (k^{2} - 1)\,(2\,k + 3)}}\, - 1 \, \right)\cdot B(k) \, < \, 0 \\[1.1em]
\!&\!\!\Longleftrightarrow\!\!&\!
- {\displaystyle \frac{4\,k^{3} + 6\,k^{2} - (c^{2} + 4)\,k - 6}{2\,(k^{2} - 1)\,(2\,k + 3)}} \cdot B(k) \, < \, 0 \\[1.1em]
\!&\!\!\Longleftrightarrow\!\!&\!
- {\displaystyle \frac{2\,(2\,k^{3} - 3)+4 \, k \,(k-1) + k (2 \, k-c^{2})}{2\,(k^{2} - 1)(2\,k + 3)}} \cdot B(k) \, < \, 0.
\end{array}
$$
As the last inequality holds for every  $k \!\geq\! 2$, the assertions of Theorem~\ref{improvement-s-2}
 follow from Proposition~\ref{Leibniz}. \hfill $\Box$


\medskip
\subsection{ Refinements of the inequalities in Statement \mbox{3}}

\smallskip
We propose the following improvement and generalization of Statement 3:

\smallskip
\begin{theorem}
\label{improvement-s-3} ~\\
$(i)$ For every $ x \!\in\! \left(0, \mbox{\small $\displaystyle\frac{\pi}{2}$}\right)$ and every $n \in N$, we have
\begin{equation}
\label{t-3} \mbox{\small $3 + \displaystyle \frac{1}{\cos x} \!\! \sum\limits_{k=2}^{2n+1}{
(-1)^{k}\, C(k)\, x^{2k}}
 \, < \,
\displaystyle 2 \, \frac{\sin{x}}{x} + \frac{\mbox{\rm tan}\,x}{x}\, < \,
3 \, + \,\displaystyle \frac{1}{\cos{x}} \! \sum\limits_{k=2}^{2n}{
(-1)^{k}\, C(k) \, \, x^{2k}}$},
\end{equation}
where $C(k)= {\displaystyle 2\,\frac{4^{k} - 3\,k - 1}{(2\,k + 1)\mathrm{!}}
} $.

\smallskip
\noindent $(ii)$ For every $\displaystyle \, x \!\in\! \left(0, \mbox{\small $\displaystyle\frac{\pi}{2}$}\right)$
and every $m \in N, \, m\geq 2$, we have the following error estimation:
$$
\begin{array}{l}
{\Big |}\,\displaystyle
2 \, \frac{\sin{x}}{x} + \frac{\mbox{\rm tan}\,x}{x}
\!-\! \left(3 + \frac{1}{\cos \, x}
\displaystyle\sum_{k=2}^{m}{(-1)^{k+1}
C(k) \, x^{2k}} \right) \,{\Big |} \, < \, C(m+1)\, \frac{x^{2m+2}}{\cos \, x}.
\end{array}
$$
\end{theorem}

\noindent
{\bf Examples.} \\
Let $x \!\in\! \left(0, \mbox{\small $\displaystyle\frac{\pi}{2}$}\right)$. \\
For $n = 1$ we get Statement~\ref{statement_3}. \\[1.0 ex]
For $n > 1 $ we have the following new results:
\begin{itemize}
\item Taking $n=2$ in (\ref{t-3}) gives:

$
2
+
\mbox{\small $\displaystyle\frac{1}{\cos x}$}
\left(
\mbox{\small $\displaystyle\frac{3}{20}$}\,x^{4}
-
\mbox{\small $\displaystyle\frac{3}{140}$}\,x^{6}
+
\mbox{\small $\displaystyle\frac{3}{2240}$}\,x^{8}
-
\mbox{\small $\displaystyle\frac{1}{19800}$}\,x^{10}\right)
 \, < \,
2 \,
\mbox{\small $\displaystyle\frac{\sin x}{x}$}
+
\mbox{\small $\displaystyle\frac{\tan x}{x}$}
\, < \, \quad \quad
< 2
+
\mbox{\small $\displaystyle\frac{1}{\cos x}$}
\left(
\mbox{\small $\displaystyle\frac{3}{20}$}\,x^{4}
-
\mbox{\small $\displaystyle\frac{3}{140}$}\,x^{6}
+
\mbox{\small $\displaystyle\frac{3}{2240}$}\,x^{8}
\right)
$

\item Taking $n=3$ in (\ref{t-3}) gives:

$
2
+
\mbox{\small $\displaystyle\frac{1}{\cos x}$}
\left(
\mbox{\small $\displaystyle\frac{3}{20}$}\,x^{4}
-
\mbox{\small $\displaystyle\frac{3}{140}$}\,x^{6}
+
\mbox{\small $\displaystyle\frac{3}{2240}$}\,x^{8}
-
\mbox{\small $\displaystyle\frac{1}{19800}$}\,x^{10}
+
\mbox{\small $\displaystyle\frac{151}{115315200}$}\,x^{12}
-
\mbox{\small $\displaystyle\frac{101}{4036032000}$}\,x^{14}
\right)
<
$

$
<
2
\mbox{\small $\displaystyle\frac{\sin x}{x}$}
+
\mbox{\small $\displaystyle\frac{\tan x}{x}$}
<
2
+
\mbox{\small $\displaystyle\frac{1}{\cos x}$}
\left(
\mbox{\small $\displaystyle\frac{3}{20}$}\,x^{4}
-
\mbox{\small $\displaystyle\frac{3}{140}$}\,x^{6}
+
\mbox{\small $\displaystyle\frac{3}{2240}$}\,x^{8}
-
\mbox{\small $\displaystyle\frac{1}{19800}$}\,x^{10}
+
\mbox{\small $\displaystyle\frac{151}{115315200}$}\,x^{12}
\right)
$

\item Taking $n=4$ in (\ref{t-3}) gives:

$
\displaystyle 2+ \frac{1}{\cos{x}}\left({\displaystyle \frac{3}{20}}\,x^{4}
-
{\displaystyle \frac{3}{140}}\,x^{6}
+
{\displaystyle \frac{3}{2240}}\,x^{8}
-
{\displaystyle \frac{1}{19800}}\,x^{10}
+
{\displaystyle \frac{151}{115315200}}\,x^{12}
-
{\displaystyle \frac{101}{4036032000}}\,x^{14}
+\right.
$

$
+\left.
{\displaystyle \frac{7279}{19760412672000}}\,x^{16}
-
{\displaystyle \frac{809}{187723920384000}}\,x^{18}
\right)
\, < \,
\displaystyle 2 \, \frac{\sin{x}}{x} + \frac{\mbox{\rm tan}\,x}{x}\, < \,
$

$
2+ {\displaystyle \frac{1}{\cos{x}}}
\left({\displaystyle \frac{3}{20}}\,x^{4}
-
{\displaystyle \frac{3}{140}}\,x^{6}
+
{\displaystyle \frac{3}{2240}}\,x^{8}
-
{\displaystyle \frac{1}{19800}}\,x^{10}
+
{\displaystyle\frac{151}{115315200}}\,x^{12}
- \right.
$

$
\left.
-{\displaystyle \frac{101}{4036032000}}\,x^{14}
+
{\displaystyle \frac{7279}{19760412672000}}\,x^{16}\right)
$
\item[etc.]
\end{itemize}

\noindent
{\bf Proof of Theorem~\ref{improvement-s-3}.}

\smallskip
Consider the function
$$
f(x) = \frac{\sin{2 x}}{x}\,+\, \frac{\sin{x}}{x} - 3 \cos{x}
$$
for $x \!\in\! \left(0, \mbox{\small $\displaystyle\frac{\pi}{2}$}\right)$.
Based on  power series expansion (\ref{sin}) we have
$$
\begin{array}{lcl}f(x)
 \!&\!\!=\!\!&\!
\displaystyle \sum _{k=0}^{\infty}
\,(-1)^{k}\,\left({\displaystyle \frac{2^{(2\,k + 1)}}{(2\,k + 1)
\mathrm{!}}} + {\displaystyle \frac{1}{(2\,k + 1)\mathrm{!}}} - {\displaystyle \frac{3}{(2\,k)\mathrm{!}}} \right)\,x^{2\,k}, \\ [1.3em]
 \!&\!\!=\!\!&\! \displaystyle \sum _{k=2}^{\infty} \,(-1)^{k}\,{\displaystyle 2 \, \frac{\,4^{k} - 3\,k - 1}{(2\,k + 1)\mathrm{!}}}\,x^{2\,k},
\end{array}
$$
for $x \!\in\! \left(0, \mbox{\small $\displaystyle\frac{\pi}{2}$}\right)$.

\smallskip
For $c=\pi/2$ we have:
$$
C(k)>0 \;\; \mbox{for} \; k\geq 2, \;\; \quad \mbox{and} \; \; \quad \displaystyle\lim_{k\rightarrow \infty }c^{2k}C(k)=0.
$$
Also,
$$
\begin{array}{l}\displaystyle
c^{2}\, C(k + 1) < C(k) \,\, \Longleftrightarrow \\[1.0 ex]
-2 \left( 12\,k^{3} \!+\! 38\,k^{2} \!+\!
{\displaystyle \frac{(6\,c^{2} \!+\! 76)\,k}{2}}
\!+\! 4\,c^{2} \!+\! 12 \right)
\,<\,
4^{k + 1}\left(2\,(k+c)(k-c) \!+\! 5\,k \!+\! 3\right).
\end{array}
$$
As the last inequality holds for $k \geq 2$, the assertions of Theorem~\ref{improvement-s-2}
immediately follows from Proposition~\ref{Leibniz}.
\hfill $\Box$


\medskip
\subsection{ Refinements of the inequalities in Statement \mbox{4}}

\smallskip
We propose the following improvement and generalization of Statement 4:

\smallskip
\begin{theorem} ~\\
\label{improvement-s-4}
$(i)$ For every $x \!\in\! \left(0, \mbox{\small $\displaystyle\frac{\pi}{2}$}\right)$
and every $n \in N$, we have
\begin{equation}
\label{t-4}
\begin{array}{c}
\displaystyle
2 + \frac{1}{\cos x}\sum_{k=2}^{2n+1} (-1)^{k} D(k) \, x^{2k}
 <
 \left(\frac{\sin{\,x}}{x}\right)^{\!2} + \frac{\tan{x}}{x}\, < \,
\displaystyle 2 + \frac{1}{\cos x}
\sum_{k=2}^{2n}(-1)^{k} D(k) \, x^{2k},
\end{array}
\end{equation}
where
$\displaystyle
D(k) = {\displaystyle \frac{1}{4}}\,{\displaystyle \frac{
 - 9 + 3^{2k + 2} - 40\,k - 32\,k^{2}}{(2\,k + 2)\mathrm{!}}}.
$

\smallskip

\noindent $(ii)$ For every $x \!\in\! \left(0, \mbox{\small $\displaystyle\frac{\pi}{2}$}\right)$
and every $m \in N, m\geq 2$, we have the following
error estimation:
$$
\begin{array}{l}
{\Big |}\displaystyle
\left(\frac{\sin{\,x}}{x}\right)^{\!\!2} \! + \frac{\tan{x}}{x} - \left(2+\frac{1}{\cos \, x}
\displaystyle\sum_{k=2}^{m}{(-1)^{k+1}
D(k) \, x^{2k}}\right){\Big |} \, < \,D(m+1)\,\frac{x^{2m+2}}{\cos \, x}.
\end{array}
$$
\end{theorem}

\noindent
{\bf Examples.} \\
Let $x \!\in\! \left(0, \mbox{\small $\displaystyle\frac{\pi}{2}$}\right)$. \\
For $n = 1$ we get Statement~\ref{statement_4}. \\[1.0 ex]
For $n > 1$ we have the following new results:
\begin{itemize}
\item Taking $n=2$ in (\ref{t-4}) gives:

$
2
+
\mbox{\small $\displaystyle\frac{1}{\cos x}$}
\left(
\mbox{\small $\displaystyle\frac{8}{45}$}\,x^{4}
-
\mbox{\small $\displaystyle\frac{4}{105}$}\,x^{6}
+
\mbox{\small $\displaystyle\frac{19}{4725}$}\,x^{8}
-
\mbox{\small $\displaystyle\frac{37}{133650}$}\,x^{10}
\right)
<
\left(
\mbox{\small $\displaystyle\frac{\sin x}{x}$}
\right)^{\!2}
+
\mbox{\small $\displaystyle\frac{\tan x}{x}$}
<
$
$
<
2
+
\mbox{\small $\displaystyle\frac{1}{\cos x}$}
\left(
\mbox{\small $\displaystyle\frac{8}{45}$}\,x^{4}
-
\mbox{\small $\displaystyle\frac{4}{105}$}\,x^{6}
+
\mbox{\small $\displaystyle\frac{19}{4725}$}\,x^{8}
\right)
$

\item Taking $n=3$ in (\ref{t-4}) gives:

$
2
+
\mbox{\small $\displaystyle\frac{1}{\cos x}$}
\left(
\mbox{\small $\displaystyle\frac{8}{45}$}\,x^{4}
-
\mbox{\small $\displaystyle\frac{4}{105}$}\,x^{6}
+
\mbox{\small $\displaystyle\frac{19}{4725}$}\,x^{8}
-
\mbox{\small $\displaystyle\frac{37}{133650}$}\,x^{10}
+
\mbox{\small $\displaystyle\frac{283}{20638800}$}\,x^{12}
-
\mbox{\small $\displaystyle\frac{3503}{6810804000}$}\,x^{14} \right)
<
$
$
<
\left(
\mbox{\small $\displaystyle\frac{\sin x}{x}$}
\right)^{\!2}
+
\mbox{\small $\displaystyle\frac{\tan x}{x}$}
\, < \,
2
+
\mbox{\small $\displaystyle\frac{1}{\cos x}$}
\left(
\mbox{\small $\displaystyle\frac{8}{45}$}\,x^{4}
-
\mbox{\small $\displaystyle\frac{4}{105}$}\,x^{6}
+
\mbox{\small $\displaystyle\frac{19}{4725}$}\,x^{8}
-
\mbox{\small $\displaystyle\frac{37}{133650}$}\,x^{10}
+
\mbox{\small $\displaystyle\frac{283}{20638800}$}\,x^{12}\right)\!
$

\item Taking $n=4$ in (\ref{t-4}) gives:

$
2
+
\mbox{\small $\displaystyle\frac{1}{\cos x}$}
\left(
\mbox{\small $\displaystyle\frac{8}{45}$}\,x^{4}
-
\mbox{\small $\displaystyle\frac{4}{105}$}\,x^{6}
+
\mbox{\small $\displaystyle\frac{19}{4725}$}\,x^{8}
-
\mbox{\small $\displaystyle\frac{37}{133650}$}\,x^{10}
+
\mbox{\small $\displaystyle\frac{283}{20638800}$}\,x^{12}
-
\mbox{\small $\displaystyle\frac{3503}{6810804000}$}\,x^{14}
+
\right.
$

$
\left.
+\mbox{\small $\displaystyle\frac{189169}{12504636144000}$}\,x^{16}
-
\mbox{\small $\displaystyle\frac{18917}{52797352608000}$}\,x^{18}
\right)
\, < \,
\left(
\mbox{\small $\displaystyle\frac{\sin x}{x}$}\right)^{\!2}
+
\mbox{\small $\displaystyle\frac{\tan x}{x}$}
\, <$

$
<
2
+
\mbox{\small $\displaystyle\frac{1}{\cos x}$}
\!\left(
\mbox{\small $\displaystyle\frac{8}{45}$}\,x^{4}
-
\mbox{\small $\displaystyle\frac{4}{105}$}\,x^{6}
+
\mbox{\small $\displaystyle\frac{19}{4725}$}\,x^{8}
-
\mbox{\small $\displaystyle\frac{37}{133650}$}\,x^{10}
+
\mbox{\small $\displaystyle\frac{283}{20638800}$}\,x^{12} -
\right.
$

$
\left.
-\mbox{\small $\displaystyle\frac{3503}{6810804000}$}\,x^{14}
+
\mbox{\small $\displaystyle\frac{189169}{12504636144000}$}\,x^{16}\right)\!
$

\item[etc.]

\end{itemize}

\medskip
\noindent
{\bf Proof of Theorem~\ref{improvement-s-4}.}

\smallskip
Consider the function
$$
f(x)
=
{\displaystyle \frac{\mathrm{cos}\,x
-
\mathrm{cos}^{3}\,x}{x^{2}}}
+
{\displaystyle \frac{\mathrm{sin}\,x}{x}}
-
2\,\mathrm{cos}\,x , \quad \mbox{for} \; x \!\in\! \left(0, \mbox{\small $\displaystyle\frac{\pi}{2}$}\right).
$$
Based on  power series expansion (\ref{cos3}) we have
$$
\begin{array}{lcl}
f(x) \!&\!\!=\!\!&\! {\displaystyle \sum _{k=0}^{\infty}}
\,(-1)^{k}\left( - {\displaystyle \frac{1}{(2\,k + 2)\mathrm{!}}
} + {\displaystyle \frac{1}{4}}\,{\displaystyle \frac{3^{(2\,
k + 2)} + 3}{(2\,k + 2)\mathrm{!}}} + {\displaystyle \frac{1}{(
2\,k + 1)\mathrm{!}}} - {\displaystyle \frac{2}{(2\,k)\mathrm{!
}}} \right)x^{2\,k} \\ [1.1em]
 \!&\!\!=\!\!&\! {\displaystyle \sum _{k=0}^{\infty}}\,(-1)^{k} \,{\displaystyle \frac{
9^{k + 1} - (32\,k^{2} + 40\,k +9)}{4 \, (2\,k + 2)\mathrm{!}}}\,x^{2\,k} \\ [1.1em]
 \!&\!\!=\!\!&\! {\displaystyle \sum _{k=2}^{\infty}}\,(-1)^{k}\,{\displaystyle \frac{
9^{k + 1} - (32\,k^{2} + 40\,k +9) }{4 \, (2\,k + 2)\mathrm{!}}}\,x^{2\,k}.
\end{array}
$$
For $ c=\pi/2$ we have:
$$
D(k)>0 \;\; \mbox{for} \; k\geq 2, \;\; \quad
\mbox{and} \; \; \quad \displaystyle \lim _{k\rightarrow \infty }c^{2k}D(k) =0.
$$
Also,

\smallskip
$
\begin{array}{l}\displaystyle
c^{2}\, D(k + 1) < D(k) ~ \Longleftrightarrow \\[0.5 ex]
{\displaystyle 32\,k^{4} + 152\,k^{3} + ( 245- 8 c^{2}
)\,k^{2} + \left({\displaystyle \frac{303}{2}} - 26\,c^{2}\right)\,k + 27
 - {\displaystyle \frac{81\,c^{2}}{4}} } \\[0.5 ex]
\hfill <
{\displaystyle \frac{1}{4}}\,{\displaystyle 9^{k+1}\,( 4 \,k^{2} + 14\,k + 12-
81\,c^{2} )}.
\end{array}
$

\medskip
\noindent
As the last inequality holds for $k \geq 2$ the assertions of Theorem~\ref{improvement-s-2}
immediately follows from Proposition~\ref{Leibniz}.
\hfill $\Box$


\medskip
\subsection{ Refinements of the inequalities in Statement \mbox{5}}

\smallskip
We prove the following generalization of Statement 5:

\begin{theorem} ~\\
\label{improvement-s-5}
$\!\!$For every $x \!\in\! \left(0, \mbox{\small $\displaystyle\frac{\pi}{2}$}\right)$ and $m \in N$, $m \ge 2$, the
following inequalities hold:
\begin{equation}
\label{t-5}
\begin{array}{c}\displaystyle
2 + \mathop \sum \limits_{k = 2}^{m - 1} \frac{{\left| {{\mbox{\boldmath $B$}_{2k}}} \right|(2k - 2){4^k}}}{{(2k)!}}{x^{2k}} +
{\left( {\frac{{2x}}{\pi}} \right)^{\!2n}}
\left( {
\frac{{{\pi^2}}}{4}
- 2 - \mathop \sum \limits_{k = 2}^{m - 1} \frac{{\left| {{\mbox{\boldmath $B$}_{2k}}} \right|(2k - 2){4^k}}}{{(2k)!}}{{\left(
{\frac{\pi}{2}} \right)}^{2k}}} \right)      \\[3.0 ex]
> \, \displaystyle {\left( {\frac{x}{{\sin x}}} \right)^2} + \frac{x}{{\tan x}}\, > \, 2 + \sum \limits_{k = 2}^m \frac{{\left| {{\mbox{\boldmath $B$}_{2k}}} \right|(2k - 2){4^k}}}{{(2k)!}}{x^{2k}}
\end{array}
\end{equation}
where $\mbox{\boldmath $B$}_{i}$ are {\sc Bernoulli}'s numbers.
\end{theorem}

\noindent
{\bf Proof of Theorem~5.}

Consider the function
$$
f(x)
=
{\displaystyle \left( {\frac{x}{{\sin\,x}}} \right)^{\!2}} + \frac{x}{{\tan\,x}},
\quad \mbox{for} \;\; x \!\in\! \left(0, \mbox{\small$\displaystyle\frac{\pi}{2}$}\right).
$$
Based on  series expansion (\ref{cosec2}) and (\ref{cotan}), we have
$$
\begin{array}{lcl}
f(x)
\!&\!\!=\!\!&\!
\displaystyle {x^2}\left( {\frac{1}{{{x^2}}} + \displaystyle \sum \limits_{k = 1}^{  \infty }
\frac{{\left| {{\mbox{\boldmath $B$}_{2k}}} \right|(2k - 1){4^k}}}{{(2k)!}}{x^{2k - 2}}} \right) + x\left( {\frac{1}{x} - \mathop
\sum \limits_{k = 1}^{  \infty } \frac{{\left| {{\mbox{\boldmath $B$}_{2k}}} \right|{4^k}}}{{(2k)!}}{x^{2k - 1}}} \right)\\[1.1em]
\!&\!\!=\!\!&\!
\displaystyle 2 + \mathop \sum \limits_{k = 1}^{  \infty } \frac{{\left| {{\mbox{\boldmath $B$}_{2k}}} \right|(2k -
1){4^k} - {4^k}}}{{(2k)!}}{x^{2k}}      \\[1.1em]
\!&\!\!=\!\!&\!
\displaystyle 2 + \mathop \sum \limits_{k = 2}^{  \infty } \frac{{\left| {{\mbox{\boldmath $B$}_{2k}}} \right|(2k -
2){4^k}}}{{(2k)!}}{x^{2k}}\\ [1.1em]
\!&\!\!>\!\!&\! \displaystyle 2 + \mathop \sum \limits_{k = 2}^m \frac{{\left| {{\mbox{\boldmath $B$}_{2k}}} \right|
(2k - 2){4^k}}}{{(2k)!}}{x^{2k}}.
\end{array}
$$
Since all coefficients are positive, by applying  Theorem WD we get the inequalities in the statement of
the theorem.
\hfill $\Box$

\medskip
\noindent
{\bf Examples.} \\
For $x \!\in\! \left( {0, \displaystyle\frac{\pi}{2} } \right)$ and
$
f(x)\!=\!{\left(\displaystyle\frac{x}{{\sin x}} \right)^{\!2}} + \displaystyle\frac{x}{{\tan x}},
$
we show the inequalities for $m \!=\! 2, 3,$ $4, 5$.
\begin{itemize}
\item For $m=2$:
\[\displaystyle
2
+
{\left( {\mbox{\small $\displaystyle\frac{2}{\pi}$}} \right)^{\!4}}\!\left(
-2
+
\mbox{\small $\displaystyle\frac{{{\pi^2}}}{4}$} \right){x^4}
> f(x) >
2
+
\mbox{\small $\displaystyle\frac{2}{{45}}$}{x^4}
\]
On the right-hand side we see the inequality from Statement 5.
\item For $m = 3$:
\[\displaystyle
2
+
\mbox{\small $\displaystyle\frac{2}{{45}}$}{x^4}
+
{\left({\mbox{\small $\displaystyle\frac{2}{\pi}$}}\right)^{\!6}}\!\left(
-2
+
\mbox{\small $\displaystyle\frac{{{\pi^2}}}{4}$} -
\mbox{\small $\displaystyle\frac{{{\pi^4}}}{{360}}$}\right){x^6}
> f(x) >
2
+ \mbox{\small $\displaystyle\frac{2}{{45}}$}{x^4}
+
\mbox{\small $\displaystyle\frac{8}{{945}}$}{x^6}.
\]
\item For $m = 4$:\\
$\begin{array}{c}\displaystyle
2
+
\mbox{\small $\displaystyle\frac{2}{{45}}$}{x^4}
+
\mbox{\small $\displaystyle\frac{8}{{945}}$}{x^6}
+
{\left(\mbox{\small $\displaystyle\frac{2}{\pi}$}\right)^{\!8}}\!\left(
-2
+
\mbox{\small $\displaystyle\frac{{{\pi^2}}}{4}$}
-
\mbox{\small $\displaystyle\frac{{{\pi^4}}}{{360}}$}
-
\mbox{\small $\displaystyle\frac{{{\pi^6}}}{{7560}}$} \right){x^8} >
\displaystyle f(x) \\ [1.3em]
\hfill >
\displaystyle
2
+
\mbox{\small $\displaystyle\frac{2}{{45}}$}{x^4}
+
\mbox{\small $\displaystyle\frac{8}{{945}}$}{x^6}
+
\mbox{\small $\displaystyle\frac{2}{{1575}}$}{x^8}.
\end{array}
$
\item For $m = 5$: \\
\[
\begin{array}{c}\displaystyle
2
\!+\!
\mbox{\small $\displaystyle\frac{2}{{45}}$}{x^4}
\!+\!
\mbox{\small $\displaystyle\frac{8}{{945}}$}{x^6}
\!+\!
\mbox{\small $\displaystyle\frac{2}{{1575}}$}{x^8}
\!+\!
{\left( \mbox{\small $\displaystyle\frac{2}{\pi}$}
\right)^{\!10}}\!\left(
\!-\!
2
\!+\!
\mbox{\small $\displaystyle\frac{{{\pi^2}}}{4}$}
\!-\!
\mbox{\small $\displaystyle\frac{{{\pi^4}}}{{360}}$}
\!-\!
\mbox{\small $\displaystyle\frac{{{\pi^6}}}{{7560}}$}
\!-\!
\mbox{\small $\displaystyle\frac{{{\pi^8}}}{{201600}}$}
\right){x^{10}} >  \\[3.0 ex]
> f(x) >
2
\!+\!
\mbox{\small $\displaystyle\frac{2}{{45}}$}{x^4}
\!+\!
\mbox{\small $\displaystyle\frac{8}{{945}}$}{x^6}
\!+\!
\mbox{\small $\displaystyle\frac{2}{{1575}}$}{x^8}
\!+\!
\mbox{\small $\displaystyle\frac{{16}}{{93555}}$}{x^{10}}.
\end{array}
\]
\end{itemize}
\noindent {\bf Remark.}\\
Let us notice that Theorem WD allows for the  approximation error to be estimated.
The difference between the right-hand side and the left-hand side of the double inequality in
Theorem~\ref{improvement-s-5} can be represented by the following function:
\[
\displaystyle
{R_n}\left( x \right)
=
\left(\!{f\left( {\frac{\pi}{2}} \right)
- 2
- \mathop \sum \limits_{k = 1}^n \frac{{\left| {{\mbox{\boldmath $B$}_{2k}}} \right|(2k - 2){4^k}}}{{(2k)!}}{{\left( {\frac{\pi}{2}} \right)}^{2k}}} \right)
\!{\left( {\frac{{2x}}{\pi}} \right)^{\!2n }}
\]
The maximum values of these functions are reached at $\mbox{\small $\displaystyle\frac{\pi}{2}$}$.
They are shown in the following table for $n = 3,4,5,6$:
$$
\begin{array}{|c|c|c|c|c|c|c|} \hline
n \!&\! 3 \!&\! 4 \!&\! 5 \!&\! 6 \\ \hline
R_{n}(x) \!&\! 0.0696524 \!&\! 0.0225863 \!&\! 0.00694574 \!&\! 0.00206173 \\ \hline
\end{array}
$$

\smallskip
\subsection{ Refinements of the inequalities in Statement \mbox{6}}

\smallskip
We propose the following generalization of Statement 6.

\begin{theorem} ~\\
\label{improvement-s-6 }
$\!\!$
For every $x \!\in\! \left( {0,\displaystyle\frac{\pi}{2}} \right)$ and $m \in N$, $m \ge 3$, the following inequality holds:
\begin{equation}
\label{t-6}
\!\!\!\!\!\!\!\!\!\!\!\!\!\!\!\!\!\!\!\!\!\!\!
\begin{array}{c}
\displaystyle
4 \!+\! \mathop \sum \limits_{k = 1}^{m - 1} \frac{3{\left| {{\mbox{\boldmath $B$}_{2k}}} \right|({2^{2k}} \!-\! 2)
\!+\! {{(-1)}^k}}}{{(2k)!}}{x^{2k}} \!+\! {{{\left(\displaystyle\frac{2x}{\pi}\right)^{\!\!2m}}}}\!\!\left(\!
{f\!\left(\frac{\pi}{2}\right) \!- 4 -\! \displaystyle \sum \limits_{k = 1}^{m - 1} \frac{3{\left|
{{\mbox{\boldmath $B$}_{2k}}}\right|({2^{2k}} \!-\! 2) \!+\! {{( - 1)}^k}}}{{(2k)!}}{{\left( {\frac{\pi}{2}} \right)}^{\!2k}}} \right) \\[1.3em]
> \displaystyle 3\frac{x}{{\sin x}} \!+\! \cos x > 4 + \mathop \sum \limits_{k = 1}^m \frac{{3\left| {{\mbox{\boldmath $B$}_{2k}}}
\right|({2^{2k}} \!-\! 2) \!+\! {{( - 1)}^k}}}{{(2k)!}}{x^{2k}}.
\end{array}
\end{equation}
where $\mbox{\boldmath $B$}_{i}$ are {\sc Bernoulli}'s numbers.
\end{theorem}
{\bf Proof of Theorem~6.}
Consider the function
$$
f(x)= \displaystyle 3\frac{x}{{\sin x}} + \cos x, \quad \mbox{for} \; x \!\in\! \left(0, \mbox{\small $\displaystyle\frac{\pi}{2}$}\right).
$$
Based on the series expansion (\ref{sin}) and (\ref{cosec}), we have
\[\displaystyle f(x)
= 4 + \displaystyle \sum \limits_{k = 1}^{ \infty } \frac{{3\left| {{\mbox{\boldmath $B$}_{2k}}} \right|({2^{2k}} - 2) + {{( -
1)}^k}}}{{(2k)!}}{x^{2k}}.\]
It is easy to verify that $3\left| {{\mbox{\boldmath $B$}_{2k}}} \right|({2^{2k}} \!-\! 2) > 1$ for $k\geq2$, and that it is equal to $1$
for $k\!=\!1$, therefore all corresponding coefficients are positive. Now, using  Theorem WD, we get the inequalities
in the statement of the theorem. \hfill $\Box$\\

\medskip
\noindent
{\bf Examples.} \\
For $x \!\in\! \left(0, \mbox{\small $\displaystyle\frac{\pi}{2}$}\right)$
and
$
f(x)=3\displaystyle\frac{x}{{\sin x}} + \cos x,
$
we show the inequalities for $m = 3, 4,$ $5, 6$.

\begin{itemize}
\item For $m = 3$:
\[
4
+
\mbox{\small $\displaystyle\frac{1}{{10}}$}{x^4}
+
{\left({\mbox{\small $\displaystyle\frac{2}{\pi}$}} \right)^{\!6}}\!\left(
-
4
+
\mbox{\small $\displaystyle\frac{{3\pi}}{2}$}
-
\mbox{\small $\displaystyle\frac{{{\pi^4}}}{{160}}$} \right){x^6}
> f(x) >
4
+
\mbox{\small $\displaystyle\frac{1}{{10}}$}{x^4}
+
\mbox{\small $\displaystyle\frac{1}{{210}}$}{x^6}.
\]
On the right-hand side we see the inequality from Statement 6.

\item For $m = 4$:
\[
4
+
\mbox{\small $\displaystyle\frac{1}{{10}}$}{x^4}
+
\mbox{\small $\displaystyle\frac{1}{{210}}$}{x^6}
+
{\left({\mbox{\small $\displaystyle\frac{2}{\pi}$}}\right)^{\!8}}\!\left(
-4
+
\mbox{\small $\displaystyle\frac{{3\pi}}{2}$}
-
\mbox{\small $\displaystyle\frac{{{\pi^4}}}{{160}}$}
-
\mbox{\small $\displaystyle\frac{{{\pi^6}}}{{13440}}$}\right){x^8}
> f(x) >
4
+
\mbox{\small $\displaystyle\frac{1}{{10}}$}{x^4}
+
\mbox{\small $\displaystyle\frac{1}{{210}}$}{x^6}
+
\mbox{\small $\displaystyle\frac{{11}}{{16800}}$}{x^8}.
\]

\item For $m = 5$:
\[
\begin{array}{c}\displaystyle
4 +
\mbox{\small $\displaystyle\frac{1}{{10}}$}{x^4}
+
\mbox{\small $\displaystyle\frac{{{1}}}{{210}}$}{x^6}
+
\mbox{\small $\displaystyle\frac{{11}}{{16800}}$}{x^8}
+
{\left({\mbox{\small $\displaystyle\frac{2}{\pi}$}}\right)^{10}}\!\left(
- 4
+
\mbox{\small $\displaystyle\frac{{3\pi}}{2}$}
-
\mbox{\small $\displaystyle\frac{{{\pi^4}}}{{160}}$}
-
\mbox{\small $\displaystyle\frac{{{\pi^6}}}{{13440}}$}
-
\mbox{\small $\displaystyle\frac{{11{\pi^8}}}{{4300800}}$} \right){x^{10}} >                                 \\[3.0 ex]
> f(x) >
\displaystyle
4
+
\mbox{\small $\displaystyle\frac{1}{{10}}$}{x^4}
+
\mbox{\small $\displaystyle\frac{1}{{210}}$}{x^6}
+
\mbox{\small $\displaystyle\frac{{11}}{{16800}}$}{x^8}
+
\mbox{\small $\displaystyle\frac{{53}}{{831600}}$}{x^{10}}.
\end{array}
\]

\item For $m = 6:$
\[\begin{array}{c}\displaystyle
4
\!+\!
\mbox{\small $\displaystyle\frac{{{x^4}}}{{10}}$}
\!+\!
\mbox{\small $\displaystyle\frac{{{x^6}}}{{210}}$}
\!+\!
\mbox{\small $\displaystyle\frac{{11{x^8}}}{{16800}}$}
\!+\!
\mbox{\small $\displaystyle\frac{{53{x^{10}}}}{{831600}}$}
\!+\!
{\left( {\mbox{\small $\displaystyle\frac{2}{\pi}$}} \right)^{\!12}}
\!\left(
-4
\!+\!
\mbox{\small $\displaystyle\frac{{3\pi}}{2}$}
\!-\!
\mbox{\small $\displaystyle\frac{{{\pi^4}}}{{160}}$}
\!-\!
\mbox{\small $\displaystyle\frac{{{\pi^6}}}{{13440}}$}
\!-\!
\mbox{\small $\displaystyle\frac{{11{\pi^8}}}{{4300800}}$}
\!-\!
\mbox{\small $\displaystyle\frac{{53{\pi^{10}}}}{{851558400}}$}
\right){x^{12}} >    \\[3.0 ex]
> f(x) >
\displaystyle
4
+
\mbox{\small $\displaystyle\frac{{{x^4}}}{{10}}$}
+
\mbox{\small $\displaystyle\frac{{{x^6}}}{{210}}$}
+
\mbox{\small $\displaystyle\frac{{11{x^8}}}{{16800}}$}
+
\mbox{\small $\displaystyle\frac{{53{x^{10}}}}{{831600}}$}
+
\mbox{\small $\displaystyle\frac{{117911{x^{12}}}}{{18162144000}}$}.
\end{array}
\]
\end{itemize}
{\bf Remark.}
The difference between the right-hand side and the left-hand side of the double inequality in Theorem 6 can be
represented by the following function:

\[
{R_n}\left( x \right)
=
\left( {f\left( {\frac{\pi}{2}} \right) - 4
+ \mathop \sum \limits_{k = 1}^n \frac{{3\left| {{\mbox{\boldmath $B$}_{2k}}} \right|({2^{2k}} - 2) + {{( -
1)}^k}}}{{(2k)!}}{{\left( {\frac{\pi}{2}} \right)}^{2k}}} \right)\!{\left( {\frac{{2x}}{\pi}} \right)^{\!2n}}.
\]
The maximum values of these functions are reached at $\mbox{\small$\displaystyle\frac{\pi}{2}$}$;
they are shown in the following table for $n = 3, 4, 5, 6$
$$
\begin{array}{|c|c|c|c|c|c|c|} \hline
n \!&\! 3 \!&\! 4 \!&\! 5 \!&\! 6 \\ \hline
R_{n}(x) \!&\! 0.0320502 \!&\! 0.00778176 \!&\! 0.00195321 \!&\! 0.000488251 \\ \hline
\end{array}
$$

\medskip
\section{Conclusion}

\smallskip
\qquad

The  idea to compare and replace functions with their
corresponding power series  to get more accurate approximations  was used in
 \cite{MRL_2017a}, \cite{MRL_2017b} and  \cite{MRL_2017c}.
 Following  the same idea, in this paper we extended the natural approach.
 We proposed and proved new inequalities which represent
refinements and generalizations of the inequalities stated in
\cite{Mortici_2011}, related to {\sc Wilker}-{\sc Cusa}-{\sc Huygens}'s inequalities.

\smallskip
Note that proofs of the new inequalities (\ref{t-1}), (\ref{t-2}), (\ref{t-3}), (\ref{t-4}), (\ref{t-5})
and (\ref{t-6}) for any fixed $n,\, m \in N$ can be obtained by substituting $x \!=\! \sin t$
for $\displaystyle t \!\in\! \left[0, \mbox{\small $\displaystyle\frac{\pi}{2}$}\right]$
and using the methods and algorithms developed in  \cite{Malesevic2016} and \cite{Lutovac2017}.
However, our approach provides proofs for the approximation of the corresponding function
by the inequality of an arbitrary degree.

The  results of the present research
 can be used to  verify  as well as to refine   a  broad
category of inequalities.
For example, {\sc Bercu} (\cite{GB_2017}, Theorem 2.4) proved  the following inequalities
$$
\!\!\!
\begin{array}{rl}
  \!&\! \left(\displaystyle\frac{x}{\sin x}\right)^2 + \displaystyle\frac{x}{\tan x}                              \\[2.0 ex]
> \!&\!
\mbox{\small $\displaystyle\frac{
\mbox{\footnotesize $11220$}\,x^{10}
\!-\!
\mbox{\footnotesize $205560$}\,x^8
\!-\!
\mbox{\footnotesize $14256000$}\,x^6
\!+\!
\mbox{\footnotesize $512179200$}\,x^4
\!-\!
\mbox{\footnotesize $3157056000$}\,x^2
\!+\!
\mbox{\footnotesize $13716864000$}}{
\mbox{\footnotesize $242$}\,x^{12}
\!-\!
\mbox{\footnotesize $8580$}\,x^{10}
\!+\!
\mbox{\footnotesize $25560$}\,x^8
\!-\!
\mbox{\footnotesize $1080000$}\,x^6
\!+\!
\mbox{\footnotesize $103680000$}\,x^4
\!-\!
\mbox{\footnotesize $1578528000$}\,x^2
\!+\!
\mbox{\footnotesize $6858432000$}}$}                                                                          \\[2.0 ex]
> \!&\! 2 + \displaystyle\frac{2}{45}x^4 > 2,
\end{array}
$$
 for every  $x \in (0, b)$ where $b=\mbox{\small $\sqrt{\displaystyle\frac{240-6\sqrt{1090}}{17}}$} = 1.5701 \ldots
< \displaystyle\frac{\pi}{2}$.

According to  Theorem~5,  taking $m=7$ in (\ref{t-5}) gives
$$\left(\displaystyle\frac{x}{\sin x}\right)^{\!2} + \displaystyle\frac{x}{\tan x} \,>\, P_{7}(x), ~~~ \mbox{for} ~ x \in \left(0, \dfrac{\pi}{2}\right)$$
where $$
P_7(x)
=
\displaystyle 2 + \mathop \sum \limits_{k = 2}^7 \frac{{\left| {{\mbox{\boldmath $B$}_{2k}}} \right|
(2k - 2){4^k}}}{{(2k)!}}{x^{2k}}.
$$
Thus {\sc Bercu}'s inequality is reduced to the following decidable problem (see \cite{Tarski_1951}, \cite{Poonen_2014}):
for every $x \in \left(0, \mbox{\small $\displaystyle\frac{\pi}{2}$}\right)$  it is true that

\vspace*{-3.5 mm}

$$
\!\!\!\!\!\!\!\!
\begin{array}{r}
  P_7(x)\! >\!
\mbox{\small $\displaystyle\frac{
\mbox{\footnotesize $11220$}\,x^{10}
\!-\!
\mbox{\footnotesize $205560$}\,x^8
\!-\!
\mbox{\footnotesize $14256000$}\,x^6
\!+\!
\mbox{\footnotesize $512179200$}\,x^4
\!-\!
\mbox{\footnotesize $3157056000$}\,x^2
\!+\!
\mbox{\footnotesize $13716864000$}}{
\mbox{\footnotesize $242$}\,x^{12}
\!-\!
\mbox{\footnotesize $8580$}\,x^{10}
\!+\!
\mbox{\footnotesize $25560$}\,x^8
\!-\!
\mbox{\footnotesize $1080000$}\,x^6
\!+\!
\mbox{\footnotesize $103680000$}\,x^4
\!-\!
\mbox{\footnotesize $1578528000$}\,x^2
\!+\!
\mbox{\footnotesize $6858432000$}}$}.
\end{array}
$$
The above inequality is a refinement of the inequality obtained by  {\sc Bercu}.
Moreover,  the above inequality   shows  that {\sc Bercu}'s inequality   holds true over an  extended interval
$\left(0, \mbox{\small $\displaystyle\frac{\pi}{2}$}\right)$.

Similarly,  Theorems 1, 2, 3, 4, 5 and 6 can be applied to other results and inequalities from \cite{GB_2016}, \cite{GB_2017}, as well as  to a broad category  of analytical inequalities.

\bigskip
\noindent \textbf{Acknowledgement.} The research of the first, second
and third authors was supported in part by the Serbian Ministry of
Education, Science and Technological Development, under Projects ON
174032, III 44006  and TR 32023 and ON 174033 respectively.

\bigskip
\noindent \textbf{Competing Interests.} The authors would like to
state that they do not have any competing interests in the subject
of this research.

\bigskip
\noindent \textbf{Author's Contributions.} All the authors
participated in every phase of the research conducted for this
paper.

\end{document}